\documentclass[12pt,fleqn]{amsart}

\usepackage{amssymb,amsmath,graphicx,enumerate}
\usepackage{color}
\usepackage{enumerate}
\usepackage[T1]{fontenc}
\usepackage{xcolor}

\newtoks\prt
\numberwithin{equation}{section}
\newtheorem{thm}{Theorem}[section]

\newtheorem{prob}[thm]{Problem}

\newtheorem{lemma}[thm]{Lemma}
\newtheorem{prop}[thm]{Proposition}
\newtheorem{cor}[thm]{Corollary}

\theoremstyle{definition}

\def\eqn#1$$#2$${\begin{equation}\label#1#2\end{equation}}

\def\A{\mathcal A}

\def\B{\mathcal B}

\def\C{\mathcal C}

\def\E{\mathcal E}

\def\V{\mathcal V}

\def\E{E^\ast}

\def\ep{\varepsilon}

\def\en{\mathbb N}

\def\er{\mathbb R}

\def \reg {\partial _{\kern1pt\text{reg}}}

\def\la{\langle}

\def\ra{\rangle}

\newcommand{\abs}[1]{\left| #1  \right|}

\newcommand{\bobo}{\beta\omega\times\beta\omega}
\newcommand{\cE}{\mathcal E}
\newcommand{\sub}{\subseteq}
\newcommand{\con}{\mathfrak c}
\newcommand{\eps}{\varepsilon}

\title{Complemented subspaces of Banach spaces $\mathcal{C}(K \times L)$}
\author[G.\ Plebanek]{Grzegorz Plebanek}
\address{Instytut Matematyczny\\ Uniwersytet Wroc\l awski\\ Pl.\ Grunwaldzki 2\\
50-384 Wroc\-\l aw\\ Poland} \email{grzegorz.plebanek@math.uni.wroc.pl}

\author[J.\ Rondo\v{s}]{Jakub Rondo\v{s}}
\address{Kurt G\"odel Research Center, Fakult\"{a}t f\"{u}r Mathematik, Universit\"{a}t Wien, Kolingasse 14-16, 1090 Wien, Austria.}
\email{jakub.rondos@gmail.com}
\author[D.\ Sobota]{Damian Sobota}
\address{Kurt G\"odel Research Center, Fakult\"{a}t f\"{u}r Mathematik, Universit\"{a}t Wien, Kolingasse 14-16, 1090 Wien, Austria.}
\email{ein.damian.sobota@gmail.com}
\thanks{The first author
was partially  supported by the NCN (National Science Centre, Poland),
under the Weave-UNISONO call in the Weave programme 2021/03/Y/ST1/00. The second author was supported by the Austrian Science Fund (FWF), Grant I 5918-N. The third author was supported by the Austrian Science Fund (FWF), Grant ESP 108-N}



\begin{document}

\begin{abstract}
We prove that, for every compact spaces $K_1,K_2$ and compact group $G$, if both $K_1$ and $K_2$ map continuously onto $G$, then the Banach space $\C(K_1 \times K_2)$ contains a complemented subspace isometric to the Banach space $\C(G)$. Consequently, $\C(K_1\times K_2)$ contains a complemented copy of $\C([0,1])$ for every non-scattered $K_1,K_2$. Also, answering a question of Alspach and Galego, we get that $\C(\beta\omega\times\beta\omega)$ contains a complemented copy of $\C([0,1]^\kappa)$ for every cardinal number $1\le\kappa\le\con$ and hence a complemented copy of $\C(K)$ for every metric compact space $K$. On the other hand, for the pointwise topology, we show that $\C_p(\beta\omega\times\beta\omega)$ contains no complemented copy of $\C_p(2^\omega)$.
\end{abstract}

\maketitle

\section{Introduction}

A Banach space $X$ is \emph{Grothendieck} if every $weak^\ast$ convergent sequence in the dual space $X^\ast$
converges weakly. Grothendieck spaces constitute an important and intensively studied class of Banach spaces, see e.g. works of Bourgain \cite{Bou83}, Diestel \cite{Die73}, Pfitzner \cite{Pfi94}, Schachermayer \cite{Sch82}, to name a few. We refer the reader to the recent extensive survey by Gonz\'alez and  Kania \cite{GK21} for basic information on Grothendieck Banach spaces.

Grothendieck \cite{Gro53} proved that for every extremally disconnected compact space $K$ the Banach space $\C(K)$, of all continuous real-valued functions on a compact space $K$ endowed with the supremum norm, is Grothendieck. This result was further generalized in various directions, in particular for the Stone spaces of Boolean algebras with some weak versions of sequential completeness, see e.g. Koszmider and Shelah \cite{KS13} for references.

It is well-known that a Banach space $\C(K)$ is Grothendieck if and only if $\C(K)$ does not contain a complemented copy of the Banach space $c_0$, see e.g.\ \cite[Proposition 4.1.2]{GK21}. Since the \v{C}ech--Stone compactification $\beta\omega$ of the set $\omega$ of natural numbers is extremally disconnected, Grothendieck's theorem yields that the space $\C(\beta\omega)$ (isometric to the Banach space $\ell_\infty$) contains no complemented copy of $c_0$. On the other hand, it follows from results of Cembranos \cite{Ce84} and Freniche \cite{Fr84} that the space $\C(K_1\times K_2)$ contains $c_0$ as a complemented subspace for every infinite compact spaces $K_1$ and $K_2$, in particular one has such a copy in the space $\C(\bobo)$. For non-separable variants of the Cembranos--Freniche theorem, requiring additional set-theoretic assumptions, see e.g. Candido and Koszmider \cite{CK16}.


Extending the results of Cembranos and Freniche, we shall prove  in Section \ref{mainproof} the following general result.

\begin{thm}
\label{main}
Let $K_1, K_2$ be compact spaces and $G$ be a compact group such that
both $K_1$ and $K_2$ map continuously onto $G$.

Then, $\C(K_1 \times K_2)$ contains a complemented subspace isometric to $\C(G)$. 
\end{thm}

Theorem \ref{main} immediately yields the following consequences.

\begin{thm}
\label{main2}
Let $K_1, K_2$ be  compact spaces such that
they both admit continuous surjections onto the cube $[0,1]^\kappa$ for some
cardinal number $\kappa\ge 1$.

Then, $\C(K_1 \times K_2)$ contains a complemented subspace isomorphic to $\C([0, 1]^\kappa)$. 
\end{thm}

\begin{proof} 
If $K_1$ and $K_2$ map continuously onto $[0,1]^\kappa$, then they also map continuously onto the product  group 
$G=T^\kappa$, where $T$ is the unit circle $[0,1)$ with addition mod 1.
Thus, by Theorem \ref{main}, $\C(K_1 \times K_2)$ contains a complemented subspace isometric to $\C(G)$. Finally, by the results of Pe{\l}czy\'nski \cite[Theorems 8.8 and 8.9]{Pe68}, $\C(G)$ is isomorphic to $\C([0,1]^{\kappa})$.
\end{proof}

Recall that a compact space $K$ is \emph{scattered} if each closed subset of $K$ contains an isolated point. By the well-known result of Pe\l czy\'{n}ski and Semadeni \cite{PS59}, $K$ is non-scattered if and only if $K$ maps continuously onto $[0,1]$ if and only if $\C(K)$ contains a copy of $\C([0,1])$. Consequently, since the product $K_1\times K_2$ of two compact spaces $K_1,K_2$ is scattered if and only if $K_1$ and $K_2$ are both scattered, the space $\C(K_1\times K_2)$ contains a copy of $\C([0,1])$ if and only if at least one of the spaces $K_1$ and $K_2$ is non-scattered. On the other hand, for $K_1, K_2$ both non-scattered, Theorem \ref{main2} yields the following supplementary result.

\begin{cor}\label{nonscattered}
If $K_1, K_2$ are non-scattered compact spaces, then $\C(K_1\times K_2)$ contains a complemented copy of $\C([0,1])$.
\end{cor} 

It follows from the aforementioned result of Cembranos and Freniche that the space $\C(\bobo)$ contains $c_0$ as a complemented subspace, even though $\C(\beta\omega)$ does not. Alspach and Galego \cite[Problem 6.4]{AG} asked if there are
other infinite-dimensional separable Banach spaces which are complemented in $\C(\bobo)$. Among non-separable ones one can immediately get the following Banach spaces: $\C(\bobo)$, $\C(\beta\omega)$, and $c_0\oplus \C(\beta\omega)$. Candido \cite[Remark 3.7]{Ca23} added the space $c_0(\C(\beta\omega))$ to this list. In the next corollary we extend significantly the collection of known separable as well as non-separable Banach spaces isomorphic to complemented subspaces of $\C(\bobo)$. 
Recall here that for any infinite discrete space $\Gamma$ its \v{C}ech--Stone compactification $\beta\Gamma$ can be continuously mapped onto every compact space of density $\le|\Gamma|$ and that the classical Hewitt--Marczewski--Pondiczery theorem implies that the density of the product $[0,1]^{\kappa}$ is $\le|\Gamma|$, for any cardinal number $\kappa\le2^{|\Gamma|}$.
 
\begin{cor}\label{cor:bobo}
For every discrete space $\Gamma$ and cardinal number $1\le\kappa\le2^{|\Gamma|}$, the space $\C(\beta\Gamma\times\beta\Gamma)$ contains a complemented subspace isomorphic to $\C([0,1]^\kappa)$. 

In particular, for any metric compact space $K$ and any cardinal number $1\le\kappa\le\con$, $\C(\bobo)$ contains a complemented copy of $\C(K^\kappa)$.
\end{cor}

The second part of the corollary follows from the facts that $\C([0,1])$ contains as a complemented subspace every separable space $\C(L)$ and that $\C(K^\kappa)$ is isomorphic to $\C([0,1]^\kappa)$ for any metric non-singleton $K$ and $\kappa\ge\omega$, see Pe\l czy\'{n}ski \cite{Pel68} and \cite[Theorem 8.8]{Pe68}. The corollary answers in particular yet another question of Alspach and Galego \cite[Problem 6.5]{AG} asking whether the space $\C([0,\omega^\omega])$ can be embedded into $\C(\bobo)$ as a complemented subspace.


In order to prove directly that $\C(\bobo)$ contains a complemented copy of $c_0$, one needs to find a bounded sequence
of signed regular Borel measures on $\bobo$ that converges in the $weak^\ast$ topology of $\C(\bobo)^\ast$ but does not converge
weakly. There are such sequences consisting of finitely supported measures, see \cite{KMSZ} and \cite{KSZ}.
However, the proof of our main result builds on another construction giving a sequence consisting of non-atomic 
probability measures on $\bobo$, cf. \cite[Theorem 9.1]{Pl23}. Our method also indicates
that the space $\C(\bobo)$ is not \emph{positively Grothendieck} in the terminology
of Koszmider and Shelah \cite{KS13}.

Corollary \ref{nonscattered} has also an immediate application to the question raised in Avil\'es \textit{et al.} \cite[Section 6.4.1]{Spaniards} asking for which compact spaces $K$ and Banach spaces $E$ the injective tensor product $\C(K)\hat{\otimes}_\eps E$ is separably injective. Recall that a Banach space $X$ is \emph{separably injective} if, for every separable Banach space $Z$ and each closed linear subspace $Y\subseteq Z$, every bounded operator $T\colon Y\to X$ extends to a bounded operator $S\colon Z\to X$. Note also that for compact spaces $K_1$ and $K_2$ the tensor product $\C(K_1)\hat{\otimes}_\eps\C(K_2)$ is isometrically isomorphic to the space $\C(K_1\times K_2)$. Clausen and Scholze \cite[Theorem 3.19]{CS22} answered the question for the special test case $\C(K)=E=\C(\beta\omega)$ by proving that the space $\C(\bobo)$ is not separably injective. Since a complemented subspace of a separably injective Banach space is itself separably injective and the space $\C([0,1])$ is not separably injective (see \cite[Chapter 2]{Spaniards}), Corollary \ref{nonscattered} yields the following result, which in particular generalizes Clausen and Scholze's one.

\begin{cor}\label{sepinj}
If $K_1, K_2$ are non-scattered compact spaces, then $\C(K_1\times K_2)$ is not separably injective.
\end{cor} 

We also briefly study the case when the function spaces are endowed with the pointwise topology. The situation appears to be completely different from that of the Banach spaces. Namely, by the theorem of Kawamura and Leiderman \cite{KL17} if, for compact spaces $K$ and $L$, the space $\C_p(L)$ is isomorphic to a complemented subspace of the space $\C_p(K)$ and $K$ has covering dimension zero, then $L$ also has covering dimension zero. In particular, $\C_p(\bobo)$ does not contain any complemented copy of $\C_p([0,1])$. We generalize the latter fact in the following way.

\begin{thm}\label{thm_cp}
    Let $K_1$ and $K_2$ be compact spaces without uncountable metrizable closed subspaces. Then, $\C_p(K_1\times K_2)$ does not contain a complemented copy of $\C_p(L)$ for any uncountable metric compact space $L$.  
\end{thm}

Theorem \ref{thm_cp} is a special case of Theorem \ref{marciszewski}. As an immediate consequence we get that $\C_p(\bobo)$ does not contain any complemented copy of $\C_p(2^\omega)$, which clearly contrasts Corollary \ref{cor:bobo}.

\section*{Acknowledgements}

The authors would like to thank Jes\'us M.F. Castillo and Witold Marciszewski for fruitful conversations and valuable comments concerning the subject of the paper.

\section{Preliminaries}

 All topological spaces considered in this paper are assumed to be Hausdorff. For a topological space $X$, $Bor(X)$ denotes the $\sigma$-field of all Borel subsets of $X$, and for a cardinal number $\kappa$ by $[X]^\kappa$, $[X]^{\le\kappa}$, and $[X]^{<\kappa}$ we denote the families of all subsets of $X$ of cardinality $\kappa$, $\le\kappa$, and $<\kappa$, respectively. For a subset $A\subseteq X$ by $\chi_A\colon X\to\{0,1\}$ we denote its characteristic function.
 
 For a compact space $K$, by $\C(K)$ we denote the Banach space of all continuous real-valued functions on $K$ endowed with the supremum norm $\|\cdot\|_\infty$. Recall that by virtue of the classical Riesz--Markov--Kakutani theorem the dual space $\C(K)^\ast$ is isometrically isomorphic to the Banach space $M(K)$ of all signed regular Borel measures of bounded variation on $K$; we will hence tacitly identify $\C(K)^\ast$ with $M(K)$, and for every $\mu\in\C(K)^\ast$ and $f\in\C(K)$ set $\mu(f)=\int_K f\,\mathrm{d}\mu$.

 Given a compact space $K$, by $P(K)$ we will denote the subspace of $\C(K)^\ast$ consisting of all probability measures. For a point $x\in K$, $\delta_x\in P(K)$ is the Dirac one-point measure at $x$.  
 By default, the space $P(K)$ is always considered with its $weak^\ast$ topology inherited from $\C(K)^\ast$.
 We refer the reader to the recent survey \cite{Pl23} on topological properties of spaces $P(K)$. 

Let $K$ be a compact space and let $\Sigma,\mathcal{A}$ be 
families of Borel subsets of $K$ such that $\Sigma \subseteq \A$.
We say that $\Sigma$ is $\bigtriangleup$\emph{-dense} in $\mathcal{A}$ with respect to some measure $\mu \in P(K)$, if for every $\ep>0$ and every $A \in \mathcal{A}$ there exists $E \in \Sigma$ such that $\mu(A \bigtriangleup E)<\ep$. 

If $g \colon K \rightarrow L$ is  a continuous surjection of compact spaces, then,
for a given measure $\mu\in P(K)$, we write $g[\mu]\in P(L)$ for the \emph{image} measure (or the \emph{push forward} measure) on $L$, that is, $g[\mu](A)=\mu(g^{-1}[A])$ for every $A\in Bor(L)$.
It is well-known that the mapping
\[P(K)\ni \mu \mapsto g[\mu]\in P(L)\]
is a continuous surjection. However, we will also need
the following result, see \cite[Theorem 2.12 and Remark 2.14]{Pl23} for the proof and further references.

\begin{prop} \label{surjection of measures}
Let $K, L$ be compact spaces and $g\colon K \rightarrow L$ be a continuous surjection. 
Then, for each $\nu \in P(L)$ there exists $\mu \in P(K)$ such that $g[\mu]=\nu$ and the family
\[\big\{g^{-1}[B] \colon B \in Bor(L)\big\}\]
 is $\bigtriangleup$-dense in $Bor(K)$ with respect to $\mu$.
\end{prop}
 
 The following lemma is a particular case of \cite[Proposition 4.1]{Pe68} characterizing continuous surjections
 admitting an averaging operator.
 
\begin{lemma} \label{projection}
Assume that $K, L$ are compact spaces such that there exist a continuous surjection $\rho \colon K \rightarrow L$ and a continuous mapping $\phi \colon L \rightarrow P(K)$ such that $\phi(y)\big(\rho^{-1}(y)\big)=1$ for each $y \in L$. Then, the mapping
\[ \C(L)\ni h\mapsto h\circ \rho\in \C(K)\]
is an isometric embedding of $\C(L)$ onto a complemented subspace of $\C(K)$.
\end{lemma}

\begin{proof}
Let the operators $S\colon \C(K)\to \C(L)$ and $T\colon \C(L)\to \C(K)$ be defined as:
\[S(f)(y)=\int_Kf(x)\,\mathrm{d}\big(\phi(y)\big)(x)\quad\text{and}\quad T(g)(x)=(g\circ\rho)(x),\]
for every $f\in \C(K)$, $y\in L$, and $g\in \C(L)$, $x\in K$. Then, $S$ is bounded and $T$ is an isometric embedding.

To finish the proof it is enough to show that $S\circ T$ is the identity on $\C(L)$. Note that for every $y\in L$ and $A\in Bor(L)$ it holds $\phi(y)(\rho^{-1}[A])=\delta_y(A)$, that is, $\rho[\phi(y)]=\delta_y$. Consequently, for every $g\in \C(L)$ and $y\in L$ we have
\[(S\circ T)(g)(y)=\int_K(g\circ\rho)(x)\,\mathrm{d}\big(\phi(y)\big)(x)=\int_Lg\,\mathrm{d}\big(\phi(y)\circ\rho^{-1}\big)=\int_Lg\,\mathrm{d}\delta_y=g(y),\]
which shows that $(S\circ T)(g)=g$ for every $g\in \C(L)$.
\end{proof}

%

We now turn to product measures. Let $K_1,K_2$ be compact spaces and let $\mu_1\in P(K_1),\mu_2\in P(K_2)$. We say that a measure $\nu\in P(K_1\times K_2)$ has the \emph{marginal distributions} $(\mu_1,\mu_2)$ if for each sets $L_1\in Bor(K_1)$ and $L_2\in Bor(K_2)$, we have
\[\nu(L_1 \times K_2)=\mu_1(L_1)\quad\text{and}\quad\nu(K_1 \times L_2)=\mu_2(L_2).\]
Recall that, in general, the $\sigma$-field $Bor(K_1)\otimes Bor(K_2)$ (generated by all Borel rectangles) may be a proper subfamily of $Bor(K_1\times K_2)$, but
for every $\mu_1\in P(K_1),\mu_2\in P(K_2)$ the product measure $\mu_1\otimes \mu_2$ has a unique extension to
a regular Borel measure on $K_1\times K_2$, so in fact $\mu_1\otimes\mu_2$ may be treated as an element of
$P(K_1\times K_2)$ (see \cite[Section 7.4]{Folland}).

We finally state the following fact concerning convergence of measures on product spaces.

\begin{lemma}
\label{rectangles}
Let $K_1, K_2$ be compact spaces and $\{ \nu_{\gamma} \}_{\gamma \in \Gamma}$ be a net in $P(K_1 \times K_2)$. Assume that all the measures $\nu_{\gamma}$ have the same marginal distributions $(\mu_1,\mu_2)\in P(K_1)\times P(K_2)$.
Let $\Sigma_1\subseteq Bor(K_1)$, $\Sigma_2\subseteq Bor(K_2)$ be such subfamilies that, for $i=1, 2$, $\Sigma_i$ is $\bigtriangleup$-dense in $Bor(K_i)$ with respect to $\mu_i$. 

Then, the following assertions are equivalent:
\begin{itemize}
    \item[(i)] the limit $\lim_{\gamma \in \Gamma} \nu_{\gamma}(A_1 \times A_2)$ exists for each $A_1 \in Bor(K_1)$ and $A_2 \in Bor(K_2)$;
    \item[(ii)] the limit $\lim_{\gamma \in \Gamma} \nu_{\gamma}(E_1 \times E_2)$ exists for each $E_1 \in \Sigma_1$ and $E_2 \in \Sigma_2$;
     \item[(iii)] the net $\{ \nu_{\gamma} \}_{\gamma \in \Gamma}$ is $weak^\ast$ convergent in $P(K_1 \times K_2)$.
\end{itemize}
In particular, a measure $\nu\in P(K_1\times K_2)$ is the $weak^\ast$ limit of $\{\nu_\gamma\}_{\gamma\in\Gamma}$ if and only if
\[\lim_{\gamma\in\Gamma}\nu_\gamma(A_1 \times A_2)=\nu(A_1\times A_2)\]
for each $A_1 \in Bor(K_1)$ and $A_2 \in Bor(K_2)$. 
\end{lemma}

\begin{proof}
Trivially, (i) implies (ii). 

For the proof of implication (ii) $\Rightarrow$ (i), assume that the limit $\lim_{\gamma \in \Gamma} \nu_{\gamma}(E_1 \times E_2)$ exists for each $E_1 \in \Sigma_1$ and $E_2 \in \Sigma_2$, and fix $\ep>0$ and an arbitrary rectangle $A_1 \times A_2 \in Bor(K_1) \otimes Bor(K_2)$. For $i=1, 2$ we find $E_i \in \Sigma_i$ such that $\mu_i(A_i \bigtriangleup E_i)<\ep/2$. Then, for each $\gamma \in \Gamma$, we have
\[|\nu_\gamma(A_1\times A_2)-\nu_\gamma(E_1\times E_2)|\le\nu_{\gamma}((A_1 \times A_2) \bigtriangleup (E_1 \times E_2))\leq\]
\[\le\nu_{\gamma}((A_1 \bigtriangleup E_1) \times K_2)+\nu_{\gamma}(K_1 \times (A_2 \bigtriangleup E_2)) =\mu_1(A_1 \bigtriangleup E_1)+\mu_2(A_2 \bigtriangleup E_2) \leq \ep,\]
hence the limit $\lim_{\gamma\in\Gamma}\nu_\gamma(A_1\times A_2)$ exists as well.

To prove that (i) implies (iii), we recall that, using the continuity and compactness, each $f \in \C(K_1 \times K_2)$ can be approximated with respect to the supremum norm by linear combinations of characteristic functions of Borel rectangles, which implies that, under the assumption of (i), also the limit $\lim_{\gamma \in \Gamma} \nu_{\gamma}(f)$ exists. Thus, by setting for each $h \in \C(K_1 \times K_2)$
\[\nu(h)=\lim_{\gamma \in \Gamma} \nu_{\gamma}(h),\]
we obtain a continuous linear functional on $\C(K_1 \times K_2)$,  which is the $weak^\ast$ limit of the net $\{ \nu_{\gamma} \}_{\gamma \in \Gamma}$. Since $\{ \nu_{\gamma} \}_{\gamma \in \Gamma}\subseteq P(K_1\times K_2)$, by the compactness of $P(K_1\times K_2)$ we have $\nu\in P(K_1 \times K_2)$.

Finally, we prove that (iii) implies (i). Thus, we assume that there exists the $weak^\ast$ limit $\nu\in P(K_1\times K_2)$ of the net $\{ \nu_{\gamma} \}_{\gamma \in \Gamma}$, and we fix Borel sets $A_1 \subseteq K_1$ and $A_2 \subseteq K_2$. Given $\ep>0$, by the regularity of the measures $\mu_1$ and $\mu_2$, for $i=1, 2$ we find an open set $U_i \subseteq K_i$ and a compact set $F_i \subseteq K_i$ such that $F_i \subseteq A_i \subseteq U_i$ and 
\[\mu_i(U_i)-\ep/2< \mu_i(A_i) < \mu_i(F_i)+\ep/2.\]
Then, similarly as above,  we have
\[\nu_{\gamma}(U_1 \times U_2)-\ep < \nu_{\gamma}(A_1 \times A_2) < \nu_{\gamma}(F_1 \times F_2)+\ep\]
for each $\gamma \in \Gamma$.
We further find a function $f \in \C(K_1 \times K_2)$, $0\le f\le 1$, such that $f$ is equal to $1$ on the set $F_1 \times F_2$, and to $0$ on the set $(K_1\times K_2)\setminus (U_1 \times U_2)$. Thus, for each $\gamma \in \Gamma$, we get
\[\nu_{\gamma}(F_1 \times F_2) \leq \nu_{\gamma}(f) \leq \nu_{\gamma}(U_1 \times U_2).\]
It follows that, for each $\gamma \in \Gamma$, we have
\[\abs{\nu_{\gamma}(f)-\nu_{\gamma}(A_1 \times A_2)} \leq \ep,\]
and hence that the limit $\lim_{\gamma\in\Gamma}\nu_\gamma(A_1\times A_2)$ exists.

\medskip

The second part of the lemma easily follows from the above proof of equivalence (i)$\Leftrightarrow$(iii). 
\end{proof}

\medskip

\section{The proof of Theorem \ref{main}}\label{mainproof}

Throughout this section, let $K_1,K_2$ be fixed compact spaces and let $G=(G,\oplus)$ be a fixed compact group. We assume that $K_1$ and $K_2$ can be continuously mapped onto $G$. Denote by $\lambda$ the unique Haar probability measure on $G$. Recall that the compactness of $G$ implies that $\lambda$ is two-sided invariant. We refer the reader to handbooks \cite[Section 11.1]{Folland} and \cite[Chapter 5]{Rudin} for general facts on Haar measures.

The proof of Theorem \ref{main} will follow in three steps.

\subsection{Measures on $G\times G$.} 

We start with the following folklore lemma. As usual, by $(L_1(G),\|\cdot\|_1)$ we denote the Banach space of all $\lambda$-integrable real-valued functions on $G$.

\begin{lemma}\label{cont}
For every $f\in L_1(G)$ the mapping $R_f\colon G\to L_1(G)$, given for every $x,y\in G$ by $R_f(y)(x)=f(x\ominus y)$, is continuous.

In particular, for every $A,B\in Bor(G)$ the mapping
\[ G\ni y\mapsto  \lambda(B\cap(A\oplus y))\]
is continuous.
\end{lemma}
\begin{proof}
Fix $f\in L_1(G)$, $\ep>0$, and $y_0\in G$. Since the space $\C(G)$ is dense in $L_1(G)$ with respect to the norm $\|\cdot\|_1$, there exists $g\in \C(G)$ such that $\|f-g\|_1<\ep/3$. It is well-known that the mapping $R_g\colon G\to \C(G)$ is continuous, see e.g. \cite[Theorem 5.13]{Rudin}. Consequently, there is a neighborhood $V$ of $y_0$ in $G$ such that for every $y\in V$ it holds
\[\big\|R_g(y)-R_g(y_0)\big\|_\infty<\ep/3.\]
For every $y\in V$ we then have
\[\big\|R_f(y)-R_f(y_0)\big\|_1=\int_G\big|f(x\ominus y)-f(x\ominus y_0)\big|\,\mathrm{d}\lambda(x)\le\]
\[\le\int_G\big|f(x\ominus y)-g(x\ominus y)\big|\,\mathrm{d}\lambda(x)+\int_G\big|g(x\ominus y)-g(x\ominus y_0)\big|\,\mathrm{d}\lambda(x)+\]
\[+\int_G\big|g(x\ominus y_0)-f(x\ominus y_0)\big|\,\mathrm{d}\lambda(x)\le2\cdot\|f-g\|_1+\big\|R_g(y)-R_g(y_0)\big\|_\infty<\ep,\]
where the penultimate inequality follows by the properties of the Haar measure, see \cite[Theorem 5.14]{Rudin}. It follows that $R_f\colon G\to L_1(G)$ is continuous for every $f\in L_1(G)$.

The second part follows from the first one by noticing that for every $A,B\in Bor(G)$ and $y\in G$ we have
\[\lambda(B\cap(A\oplus y))=\int_G\chi_{A\oplus y}(x)\cdot\chi_B(x)\,\mathrm{d}\lambda(x)=\int_G\chi_{A}(x\ominus y)\cdot\chi_B(x)\,\mathrm{d}\lambda(x)=\]
\[=\int_G R_{\chi_A}(y)(x)\cdot\chi_B(x)\,\mathrm{d}\lambda(x),\]
and that the functional
\[L_1(G)\ni f\mapsto\int_Gf\cdot\chi_B\,\mathrm{d}\lambda\]
is continuous.
\end{proof}

Next, we denote by $\mathcal{B}$ the family of all finite partitions of $G$ into Borel sets. For $\cE_1, \cE_2 \in \mathcal{B}$ we write $\cE_1\prec \cE_2$ if $\cE_2$ is finer
than $\cE_1$. Note that $(\mathcal{B},\prec)$ is an upward directed set. Further, for $\cE \in \mathcal{B}$, let $\cE^{\ast}$ stand for the family of all sets from $\cE$ that are of non-zero measure $\lambda$.

We consider the product measure $\lambda_2=\lambda\otimes\lambda$
on $G\times G$. As mentioned above, after taking the completion, $\lambda_2\in P(G\times G)$. We further consider some additional probability measures on $G\times G$.
Namely, for a given partition $\cE \in \mathcal{B}$ and an element $y\in G$, we
define
\[ \nu_{y,\cE}(D)=\sum_{E\in\cE^{\ast}}\frac{\lambda_2 \big(D\cap (E\times (E\oplus y))\big)}{\lambda(E)}\]
for every $D\in Bor(G\times G)$.

\begin{lemma}\label{proof:1}
Given $\cE\in\mathcal{B}$ and $y\in G$, we have $\nu_{y,\cE}\in P(G\times G)$ and
\[\nu_{y,\cE}(B\times G)=\nu_{y,\cE}(G\times B)=\lambda(B)\]
for every Borel set $B\sub G$.
\end{lemma}

\begin{proof}
Fix $\cE\in\mathcal{B}$ and $y\in G$. Let $B\in Bor(G)$. Since, for $E \in \cE$,
\[\lambda_2\big( (B\times G)\cap (E\times (E\oplus y))\big)= 
\lambda_2\big((B\cap E)\times (E\oplus y)\big)=\]
\[= \lambda(B\cap E)\cdot \lambda (E\oplus y)=
\lambda(B\cap E)\cdot \lambda (E),\]
we get 
\[ \nu_{y,\cE}(B\times G)=\sum_{E\in\cE^{\ast}} \lambda(B\cap E)=\sum_{E\in\cE} \lambda(B\cap E)=\lambda(B),\] 
as $\cE$ is a partition of $G$.
Likewise,
\[ \nu_{y,\cE}(G\times B)=\sum_{E\in\cE} \lambda(B\cap (E\oplus y))=\lambda(B).\]

The fact that $\nu_{y,\cE}\in P(G\times G)$ follows from the above equalities applied to $B=G$.
\end{proof}

Given $y\in G$, we write $s_y\colon G\to G\times G$ for
the mapping $s_y(x)=(x,x\oplus y)$. Let $\nu_y=s_y[\lambda]$
be the image measure of $\lambda$, that is, $\nu_y$  is the measure  $\lambda$ put
on the `shifted' diagonal 
\[ L_y=\big\{(x,x\oplus y)\colon  x\in G\big\}.\]

\begin{lemma}\label{proof:2}
Fix $y\in G$. For every Borel rectangle $A\times B\sub G\times G$ we have
\[ \nu_y(A\times B)=\lim_{\cE\in\mathcal{B}}\nu_{y,\cE}(A\times B).\]
Consequently, it holds
\[\nu_y=weak^\ast{\mathchar"2D}\lim_{\cE \in \mathcal{B}}\nu_{y,\cE},\] 
that is, the
net $\{\nu_{y,\cE}\}_{\cE\in\mathcal{B}}$ converges to $\nu_y$ in the $weak^\ast$ topology.
\end{lemma}

\begin{proof}
Consider Borel sets $A,B\sub G$ and a
partition $\cE \in \mathcal{B}$ such that $A$ is the union of
some elements of $\cE$. Then,
\[ \nu_{y,\cE}(A\times B)=\sum_{\substack{E\in\cE^{\ast}\\E\subseteq A}} 
\frac{\lambda_2\big((A\times B)\cap (E\times (E\oplus y))\big)}{\lambda(E)}=\]
\[ =\sum_{\substack{E\in\cE^{\ast}\\E\subseteq A}} 
\frac{\lambda_2\big (E\times (B \cap  (E\oplus y))\big)}{\lambda(E)}
= \sum_{\substack{E\in\cE^{\ast}\\E\subseteq A}}\lambda(B\cap (E\oplus y))=
\lambda(B\cap(A\oplus y)).
\]
On the other hand,
\[\tag{$*$}\nu_y(A\times B)=\lambda\big(\{x\in A\colon  x\oplus y\in B\}\big)=\lambda(A\cap (B\ominus y))=\]
\[=\lambda\big((A\cap(B\ominus y))\oplus y\big)=\lambda(B\cap(A\oplus y)),\]
so we get $\nu_{y,\cE}(A\times B)=\nu_y(A\times B)$ if
the cover $\cE$ is fine enough.

The convergence of the net $\{\nu_{y,\cE}\}_{\cE \in \mathcal{B}}$ to $\nu_y$ in the $weak^\ast$ topology now follows by Lemma \ref{rectangles}.
\end{proof}

\subsection{Measures on $K_1\times K_2$.}
We fix continuous surjections
\[f_i\colon K_i\to G,\ i=1,2,\]
and set 
\[f=f_1\times f_2\colon  K_1\times K_2\to G\times G,\]
so that $f$ is a continuous surjection, too. 
For $i=1,2$ we apply Lemma \ref{surjection of measures} to find a measure $\mu_i\in P(K_i)$ such that 
\begin{enumerate}[(i)]
\item $f_i[\mu_i]=\lambda$,
\item the family $\Sigma_i=\big\{f_i^{-1}[A]\colon  A\in Bor(G)\big\}$ is
$\bigtriangleup$-dense in $Bor(K_i)$ with respect to $\mu_i$.
\end{enumerate}

Now, we define certain measures on $K_1 \times K_2$ in an analogous way, as we did in the previous section on $G \times G$.
For a given partition $\cE \in \mathcal{B}$ and an element $y\in G$, we
define
\[ \mu_{y,\cE}(D)=\sum_{E\in\cE^{\ast}}\frac{\mu_1\otimes\mu_2 \Big( D\cap \big( f_1^{-1}[E]\times f_2^{-1}[E\oplus y]\big)\Big)}{\lambda(E)}\]
for each $D\in Bor(K_1 \times K_2)$. 
The proof of the next lemma is similar to the one of Lemma \ref{proof:1}.

\begin{lemma}\label{proof:3}
Fix $\cE \in \mathcal{B}$ and $y\in G$. Then, $\mu_{y,\cE}\in P(K_1\times K_2)$, $f[\mu_{y,\cE}]=\nu_{y,\cE}$, and
\[\mu_{y,\cE}(A\times K_2)=\mu_1(A), \quad \mu_{y,\cE}(K_1\times B)=\mu_2(B)\]
for every Borel sets $A\sub K_1$, $B\sub K_2$.
\end{lemma}

\begin{proof}
We start by proving the marginal equalities. We may assume that $\mu_1(A)>0$, otherwise the required equality holds trivially. If $A\sub f_1^{-1}[E]$ for some $E\in\cE$, then $E\in\cE^{\ast}$, and thus 
\[ \mu_{y,\cE}(A\times K_2)=\frac{\mu_1\otimes\mu_2\big((A\cap f_1^{-1}[E])\times(K_2\cap f_2^{-1}[E\oplus y])\big)}{\lambda(E)}=\]
\[= \frac{\mu_1(A)\cdot \mu_2\big( f_2^{-1}[E\oplus y]\big)}{\lambda(E)}= \frac{\mu_1(A)\cdot f_2[\mu_2](E\oplus y)}{\lambda(E)}=\mu_1(A).\]

Thus, for a general Borel set $A\sub K_1$, we get
\[\mu_{y,\cE}(A\times K_2)=\sum_{E \in \cE^{\ast}} \mu_{y,\cE}\big((A \cap f_1^{-1}[E])\times K_2\big)=\sum_{E \in \cE^{\ast}} \mu_1\big(A \cap f_1^{-1}[E]\big)=\mu_1(A).\]
The second formula can be checked in a similar manner.

Setting $A=K_1$, we get $\mu_{y,\cE}(K_1\times K_2)=\mu_1(K_1)=1$, so $\mu_{y,\cE}\in P(K_1\times K_2)$. Finally, to prove that $f[\mu_{y,\cE}]=\nu_{y,\cE}$, it is enough to show that for every $D\in Bor(G\times G)$ and $E\in\cE^{\ast}$ we have
\[\mu_1\otimes\mu_2 \Big(f^{-1}[D]\cap \big( f_1^{-1}[E]\times f_2^{-1}[E\oplus y]\big)\Big)=\lambda_2 \big(D\cap (E\times (E\oplus y))\big).\]
So, fix $D$ and $E$ as above. Then, we have the following standard computations, involving Fubini's theorem in the obvious two places:
\[\mu_1\otimes\mu_2 \Big(f^{-1}[D]\cap \big( f_1^{-1}[E]\times f_2^{-1}[E\oplus y]\big)\Big)=\]
\[=f[\mu_1\otimes\mu_2]\big(D\cap(E\times(E\oplus y))\big)=\int_{E\times(E\oplus y)}\chi_D\,\mathrm{d}f[\mu_1\otimes\mu_2]=\]
\[=\int_{f_1^{-1}[E]\times f_2^{-1}[E\oplus y]}\chi_D\big(f_1(x_1),f_2(x_2)\big)\,\mathrm{d}(\mu_1\otimes\mu_2)(x_1,x_2)=\]
\[=\int_{f_1^{-1}[E]}\int_{f_2^{-1}[E\oplus y]}\chi_D\big(f_1(x_1),f_2(x_2)\big)\,\mathrm{d}\mu_2(x_2)\,\mathrm{d}\mu_1(x_1)=\]
\[=\int_{E}\int_{E\oplus y}\chi_D(y_1,y_2)\,\mathrm{d}\big(f_2[\mu_2]\big)(y_2)\,\mathrm{d}\big(f_1[\mu_1]\big)(y_1)=\]
\[=\int_{E}\int_{E\oplus y}\chi_D(y_1,y_2)\,\mathrm{d}\lambda(y_2)\,\mathrm{d}\lambda(y_1)=\int_{E\times(E\oplus y)}\chi_D\,\mathrm{d}\lambda_2=\]
\[=\lambda_2\big(D\cap(E\times(E\oplus y))\big).\]
\end{proof}

\begin{lemma}\label{proof:4}
For every $y\in G$, the net $\{ \mu_{y,\cE} \}_{\cE \in \B}$ converges to some $\mu_y\in P(K_1\times K_2)$
such that $f[\mu_y]=\nu_y$.
 \end{lemma}

\begin{proof}
Fix $y\in G$. By Lemmas \ref{proof:2} and \ref{proof:3} all the measures $\mu_{y,\cE}$ have the same marginal distributions $(\mu_1,\mu_2)$ and for every $A,B\in Bor(G)$ we have
\[\lim_{\cE\in\mathcal{B}}\mu_{y,\cE}\big(f_1^{-1}[A]\times f_2^{-1}[B]\big)=\lim_{\cE\in\mathcal{B}}f[\mu_{y,\cE}](A\times B)=\lim_{\cE\in\mathcal{B}}\nu_{y,\cE}[A\times B]=\nu_y(A\times B).\]
Consequently, since the families $\Sigma_1$ and $\Sigma_2$ are $\bigtriangleup$-dense in $Bor(K_1)$ and $Bor(K_2)$, respectively, Lemma \ref{rectangles} implies that the net $\{\mu_{y,\cE} \}_{\cE \in \B}$ converges to some $\mu_y\in P(K_1\times K_2)$ in the $weak^\ast$ topology. By the continuity of the push forward operation we have $f[\mu_y]=\nu_y$.
\end{proof} 

\begin{lemma}\label{proof:5}
The mapping $\phi\colon G\to P(K_1\times K_2)$ given by the assignment
\[ G\ni y\mapsto \mu_y\in P(K_1\times K_2),\]
where each $\mu_y$ is as in Lemma \ref{proof:4}, is continuous.
\end{lemma}

\begin{proof}
Let $\{y_\gamma\}_{\gamma\in\Gamma}\subseteq G$ be a net converging to some $y\in G$. We will show that the net $\{\mu_{y_\gamma}\}_{\gamma\in\Gamma}$ converges to $\mu_y$ in the $weak^\ast$ topology. Since the families $\Sigma_1$ and $\Sigma_2$ are $\bigtriangleup$-dense in $Bor(K_1)$ and $Bor(K_2)$, respectively, by Lemma \ref{rectangles} it is enough to check that for every $A,B\in Bor(G)$ we have
\[\lim_{\gamma\in\Gamma}\mu_{y_\gamma}\big(f_1^{-1}[A]\times f_2^{-1}[B]\big)=\mu_y\big(f_1^{-1}[A]\times f_2^{-1}[B]\big).\]
So, fix $A,B\in Bor(G)$ and $\ep>0$. Similarly as in ($*$) in the proof of Lemma \ref{proof:2}, for every $y'\in G$, we have
\[ \mu_{y'}\big(f_1^{-1}[A]\times f_2^{-1}[B]\big)=\nu_{y'}(A\times B)=\lambda(B\cap(A\oplus y')),\]
so by Lemma \ref{cont} there is an open neighborhood $V$ of $y$ in $G$ such that
\[\big|\mu_{y}\big(f_1^{-1}[A]\times f_2^{-1}[B]\big)-\mu_{y'}\big(f_1^{-1}[A]\times f_2^{-1}[B]\big)\big|<\ep\]
for every $y'\in V$. Consequently, there is $\gamma_0\in\Gamma$ such that for every $\gamma\ge\gamma_0$, $\gamma\in\Gamma$, we have $y_\gamma\in V$  and hence
\[\big|\mu_{y}\big(f_1^{-1}[A]\times f_2^{-1}[B]\big)-\mu_{y_\gamma}\big(f_1^{-1}[A]\times f_2^{-1}[B]\big)\big|<\ep,\]
which finishes the proof of the lemma.
\end{proof}

\subsection{The final step.} \label{finalstep}
We consider the mapping $\rho\colon K_1\times K_2\to G$ given for every $t_1\in K_1,t_2\in K_2$ by the formula
\[\rho(t_1,t_2)=\ominus f_2(t_2)\oplus f_1(t_1).\]
It is immediate that $\rho$ is a continuous surjection.
Then, we have the isometric embedding
\[ \C(G)\ni h\mapsto h\circ \rho\in \C(K_1\times K_2).\]
Moreover, for every $y\in G$, it holds $\rho^{-1}(y)=f^{-1}[L_y]$, and so for the mapping $\phi$ from Lemma \ref{proof:5} we have
\[\phi(y)\big(\rho^{-1}(y)\big)=\mu_y\big(f^{-1}[L_y]\big)=f[\mu_y](L_y)=\nu_y(L_y)=1.\]
Consequently, by Lemma \ref{projection}, 
the above embedding is onto a complemented subspace of $\C(K_1\times K_2)$. The proof of Theorem \ref{main} is thus finished.

\section{Further results, remarks, and questions}

We start the final section of the paper by asking the following question, related to Corollary \ref{cor:bobo}. Recall here that, for any separable compact space $K$, the space $\bobo$ maps continuously onto $K$ and hence $\C(K)$ embeds isometrically into $\C(\bobo)$. We are however not aware of any example of a separable compact space $K$ such that $\C(\bobo)$ does not contain a complemented copy of $\C(K)$.

\begin{prob}\label{problem_example}
Does $\C(\bobo)$ contain a complemented isomorphic copy of $\C(K)$ for every separable compact space $K$?
\end{prob}

A natural and more general variant of Problem \ref{problem_example} could ask for a characterization of those separable compact spaces $K$ for which $\C(K)$ is not isomorphic to a complemented subspace of $\C(\bobo)$.

\subsection{Slight strengthening of Theorem \ref{main2}}

Examining briefly the proof of Theorem \ref{main2}, 
we can easily weaken its assumption. Consider a compact space $K$ containing the product
$K_1\times K_2$ of two compact spaces $K_1,K_2$ and suppose that, for $i=1,2$, a function $h_i\colon K_i\to [0,1]^\kappa$ is a continuous surjection.
Then, the product mapping \[ h=h_1\times h_2\colon  K_1\times K_2\to [0,1]^\kappa \times [0,1]^\kappa\]
 admits a
continuous extension to $\widetilde{h}\colon K\to [0,1]^\kappa \times [0,1]^\kappa$ (by
the Tietze extension theorem applied coordinatewise). Let $\sigma$ be a surjection from
$[0,1]^\kappa$ onto a suitable compact group $G$; then, using the notation from the previous section,
we have $f_i=\sigma\circ h_i$, $i=1,2$. Now, the mapping $\widetilde{\rho}\colon  K\to G$ given for every $x\in K$ by
\[\widetilde{\rho}(x)=\ominus \sigma\big(\pi_2(\widetilde{h}(x))\big) \oplus \sigma\big(\pi_1(\widetilde{h}(x))\big),\]
where both $\pi_i$'s are the natural projections, is an extension of $\rho$ from Section \ref{finalstep}, so also a surjection. Hence, $\widetilde{\rho}^{-1}(y)\supseteq \rho^{-1}(y)$ for every $y\in G$, and
we can use Lemma \ref{projection} again to conclude the following result.

\begin{cor}\label{general}
Let $K_1, K_2$ be  compact spaces such that
they admit continuous surjections onto $[0,1]^\kappa$ for some
cardinal number $\kappa\ge 1$.

If a compact space $K$ contains a subspace homeomorphic to $K_1\times K_2$,
then $\C(K)$ contains a complemented subspace isomorphic to $\C([0, 1]^\kappa)$. 
\end{cor}

Note that Corollary \ref{general} applies in particular to the compact space $K=B_{\C(K_1\times K_2)^\ast}$, the dual unit ball of a Banach space $\C(K_1\times K_2)$ endowed with the $weak^\ast$ topology.

\subsection{Zero-dimensional compact spaces}

\newcommand{\clop}{\rm clop}
\newcommand{\fA}{\mathfrak A}
\newcommand{\fB}{\mathfrak B}

If $K$ is a zero-dimensional compact space, then we write $\clop(K)$ for the Boolean algebra of clopen subsets of $K$.
In such a case we can reformulate Corollary \ref{general} using the classical
duality between Boolean algebras and their Stone spaces, see Koppelberg \cite{Ko88}.

\begin{cor}\label{zerodim}
If $K$ is a zero-dimensional compact space and $\clop(K)$
admits a homomorphic surjection onto the free product $\fA\otimes\fB$  of two non-atomic Boolean algebras $\fA$ and $\fB$, then
$\C(K)$ contains a complemented copy of $\C([0,1])$.
\end{cor}

Consider the double arrow space $S$, that  is,
\[ S=\big((0,1]\times\{0\}\big)\cup \big( [0,1)\times \{1\}\big)\]
and $S$ is equipped with the topology defined by the lexicographic order.
Recall that the natural projection $\pi\colon S\to [0,1]$ is continuous but the resulting copy $X$ of 
the space $\C([0,1])$ inside $\C(S)$ is not complemented. In fact, Kalenda and Kubi\'s \cite{KK12} proved that
$X$ is contained in no complemented separable superspace. 

On the other  hand,
Marciszewski \cite{Ma08} proved that $\C(S)$ does contain a complemented copy of $\C([0,1])$.  
We can conclude that the implication of Corollary \ref{zerodim} cannot be reversed,
since the space $S$ does not contain any product of two non-scattered compact spaces. Indeed, every perfect subset of $S$
is not metrizable, so if we suppose that a product $K_1\times K_2$ is homeomorphic to a subspace of $S$,
then one of those factors must be non-metrizable, too. Recall that  the algebra $\clop (S)$ does not allow a homomorphism onto
$\fA\otimes\fB$ where $\fA$ is uncountable---this may be derived
using the concept of minimally generated Boolean algebras, see Koppelberg \cite{Ko88},
or the notion of a free dimension introduced in \cite{MP19}.




\subsection{Spaces $\C(K_1\times K_2)$ endowed with the pointwise topology\label{sec_pt}}

Given a compact space $K$, by $\C_p(K)$ we denote the space $\C(K)$ endowed with the pointwise topology. Similarly, by $(c_0)_p$ we mean the set $\{x\in\er^\omega\colon x(n)\to0\}$ endowed with the product topology inherited from $\mathbb{R}^\omega$. We say that $\C_p(K)$ contains a \emph{complemented} copy of a topological vector space $X$ if there exist closed linear subspaces $E,F$ of $\C_p(K)$ such that $\C_p(K)=E\oplus F$ (in the algebraic sense), the canonical projections are continuous, and $E$ is linearly isomorphic to $X$. 

As mentioned in the introduction, it was shown in \cite{KSZ} that $\bobo$ admits a sequence of finitely supported measures which is $weak^\ast$ convergent to $0$ but not weakly convergent. This, together with the theorem of Banakh, K\k{a}kol, and \'{S}liwa \cite{BKS}, yields that $\C_p(\bobo)$ contains a complemented copy of $(c_0)_p$. Consequently, by virtue of the closed graph theorem, one obtains a complemented copy of the Banach space $c_0$ in $\C(\bobo)$. In the light of Corollary \ref{cor:bobo}, one can wonder whether the existence of a complemented copy of $\C([0,1])$ in $\C(\bobo)$ can also be achieved using finitely supported measures, and ask the following.

\begin{prob}\label{cp}
  Given a  metrizable compact space $L$, does the space $\C_p(\bobo)$ contain a complemented copy of $\C_p(L)$? 
\end{prob}

The answer is negative for $L=[0,1]$ by the aforementioned result due to Kawamura and Leiderman \cite{KL17} which asserts, in particular, that $\C_p(\bobo)$ does not contain a complemented copy of $\C_p(L)$ for any compact space $L$ of positive covering dimension.
It turns out that the answer to Problem \ref{cp} is negative whenever $L$ is  an uncountable metric compact space, see Theorem \ref{marciszewski}. In particular, $\C_p(\bobo)$ contains no complemented copy of the space $\C_p(2^\omega)$. The presented argument was suggested to the authors by Witold Marciszewski.

\medskip

\newcommand{\supp}{\protect{\rm supp}}
\newcommand{\vf}{\varphi}

Let $K$ and $L$ be compact spaces and let $T\colon\C_p(K)\to \C_p(L)$ be a continuous linear surjection.  Given  $y\in L$, the formula
\[\C(K)\ni g\mapsto (Tg)(y)\]
defines a continuous linear functional on the space $\C_p(K)$, which can be represented as a finite linear combination of point-evaluations,
in other words, as a measure concentrated on a finite subset $\supp_T(y)$ of $K$ (see \cite[Sections 1.3--1.4]{BG92}). We refer the reader to Baars and de Groot \cite{BG92} or Marciszewski \cite{Ma97} for basic properties of the finite-valued map $L\ni y\mapsto \supp_T(y)\in [K]^{<\omega}$. 

To discuss the mapping $\supp_T$, one can equip $[K]^{<\omega}$ with the \emph{Vietoris topology}, which, recall, is defined as follows. For any open set $U\sub K$ write
\[ U^+=\{F\in [K]^{<\omega}\colon F\sub U\},\quad U^-=\{F\in [K]^{<\omega}\colon F\cap U\neq\emptyset\}.\]
Then, sets of the form $U^+$ and $U^-$ form the subbase of the Vietoris topology. Recall also that a function $\vf\colon L\to[K]^{<\omega}$ is \emph{lower semicontinuous} if the set $\{y\in L\colon \vf(y)\in U^-\}$ is open in $L$ for every open $U\sub K$, and is \emph{upper semicontinuous} if the set $\{y\in L\colon \vf(y)\in U^+\}$ is open in $L$ for every open $U\sub K$. It is easy to see that $\supp_T$ is lower semicontinuous (see \cite[Proposition 1.4.4]{BG92}).

In the lemma given below we gather some standard information on finite-valued maps.

\begin{lemma}\label{usc}
Let $K$ and $L$ be compact spaces and let $\vf\colon L\to [K]^{<\omega}$ be a lower semicontinuous mapping for which
 there is $n\in\omega$ such that $\vf(y)\in [K]^n$ for every $y\in Y$. Then, $\vf$ is upper semicontinuous and so continuous.

Moreover, for the space $K_0=\bigcup \{\vf(y)\colon y\in L\}$ we have:
\begin{enumerate}[(i)]
\item $K_0$ is compact;
\item $K_0$ is metrizable whenever $L$ is metrizable.
\end{enumerate}
\end{lemma}

\begin{proof}
It is easy to check the upper semicontinuity of $\vf$. Clause $(i)$ is also standard, see e.g.\ \cite[Proposition 3.1]{Pl14}.
To check $(ii)$, assume that $L$ is metrizable. Write $S=\vf[L]$ and
\[ \la U_1,\ldots, U_n\ra=\{s\in S\colon s\cap U_i\neq\emptyset\mbox{ for every } i=1,\ldots,n\},\]
where $U_1,\ldots,U_n\sub K$ are open and pairwise disjoint.

As $S$ is an image of the compact metrizable space $L$ under the continuous mapping $\vf$,
$S$ has a countable base. It follows that we can fix a countable base $\V$ of $S$ consisting of
sets of the form $\la U_1,\ldots, U_n\ra $.

Consider now the mapping $\theta\colon K^n\to [K]^{\le n}$ given for every $(x_1,\ldots, x_n)\in K^n$ by 
\[ \theta(x_1,\ldots, x_n)=\{x_1,\ldots, x_n\}.\]
Clearly, $\theta $ is continuous. Hence,  $X=\theta^{-1}[S]$ is a compact subset of $K^n$.
Note that the family of all sets  $U_1\times\ldots\times U_n$, where $U_i$'s (possibly after some permutation) satisfy  $\la U_1,\ldots, U_n\ra\in\V$ 
is a countable base of $X$, so $X$ is metrizable. Finally, $K_0$ is a continuous image of $X$ (under the projection onto the first coordinate) and thus $K_0$ is metrizable, too.
%
%
(Alternatively, the fact that $X$ is metrizable follows also from the observation that every point of $X$ has a neighborhood on which $\theta$ is a homeomorphism  and that a locally metrizable compact space is metrizable.)
\end{proof}

The negative answer to Problem \ref{cp} mentioned above follows immediately from the
following more general result. 

\begin{thm}[Marciszewski] \label{marciszewski}
Let $K$ and $L$ be compact spaces. If $L$ contains a copy of the Cantor space $2^\omega$ and there is a continuous linear
surjection $T\colon\C_p(K)\to \C_p(L)$, then $K$ contains a copy of $2^\omega$ as well.
\end{thm}

\begin{proof}
We may think of the Cantor set $2^\omega$ as of a subspace of $L$ and consider the
mapping $2^\omega\ni y\mapsto \supp_T(y)\in [K]^{<\omega}$.

Note first that, for every $n\in\omega$, the set $L_n=\{y\in 2^\omega\colon |\supp_T(y)|\le n\}$ is closed. 
As $2^\omega=\bigcup_{n\in\omega} L_n$, there is the least $n_0$ such that $L_{n_0}$ has non-empty interior.
Clearly, $\supp_T(y)$ is non-empty for every $y\in 2^\omega$, so $n_0>0$. Hence,
the size of $\supp_T(y)$ is constant for every $y\in L'=L_{n_0}\setminus L_{n_0-1}$.
As $L'$ again contains a copy of the Cantor set, we reduce the setting to the case when
$\supp_T$ maps $2^\omega$ into $[L]^n$ for some fixed $n\in\omega$.

Now, we know from Lemma \ref{usc}  that $\supp_T$ is continuous with respect to the Vietoris topology and
 the set
\[ K_0=\bigcup\big\{\supp_T(y)\colon y\in 2^\omega\big\} \]
is compact and metrizable. Note also that $K_0$ is uncountable because
$2^\omega$ is uncountable (given a finite set $F\sub K$, we can have $\supp_T(y)=F$ only for finitely many $y$'s).
Hence, $K_0$ contains a copy of the Cantor set $2^\omega$, and we are done.
\end{proof}

Let us finally note that a similar argument as above, together with \cite[Lemma 3.2]{Pl14}, shows
that $\C_p(\bobo)$ does not admit a continuous linear surjection onto $\C_p(L)$ whenever
$L$ is a  compact space without isolated points and which is additionally sequentially compact.

\end{document}